\documentclass[12pt]{amsart}
\usepackage[utf8]{inputenc}

\usepackage[english]{babel}

\usepackage{amsmath,amsfonts,amssymb,amsthm}

\usepackage{amscd}

\usepackage[a4paper,nomarginpar]{geometry}

\usepackage{epsfig,graphicx,color}

\usepackage[all]{xy}

\usepackage[utf8]{inputenc}

\usepackage[english]{babel}

\usepackage{amsmath,amsfonts,amssymb,amsthm}

\usepackage{amscd}

\usepackage[a4paper,nomarginpar]{geometry}

\usepackage{epsfig,graphicx,color}

\usepackage[all]{xy}

\newtheorem{theorem}{Theorem}
\newtheorem{lemma}[theorem]{Lemma}
\newtheorem{proposition}[theorem]{Proposition}
\newtheorem{corollary}[theorem]{Corollary}

\theoremstyle{definition}
\newtheorem{definition}[theorem]{Definition}

\theoremstyle{remark}
\newtheorem{remark}[theorem]{Remark}

\newcommand{\N}{\mathbb{N}}  % set of natural numbers
\newcommand{\Z}{\mathbb{Z}}  % set of integer numbers
\newcommand{\eps}{\varepsilon} % abbreviation for epsilon
\begin{document}

\title[Shadowing and  Stability in $p$-adic dynamics]
{Shadowing and Stability in $p$-adic dynamics}

% Only \author and \address are required; other information is optional.
% Remove any unused author tags.

\author{Jéfferson Bastos, Danilo Caprio and Ali Messaoudi}
\address{Departamento de Matem\'atica, Universidade Estadual Paulista, Rua Crist\'ov\~ao Colombo, 2265, Jardim Nazareth, S\~ao Jos\'e do
    Rio Preto, SP, 15054-000, Brasil}
\email{jefferson.bastos@unesp.br, danilo.caprio@unesp.br, ali.messaoudi@unesp.br}

%\author{Ali Messaoudi}
%\address{Departamento de Matem\'atica, Universidade Estadual Paulista,
 %   Rua Crist\'ov\~ao Colombo, 2265, Jardim Nazareth, S\~ao Jos\'e do
  %  Rio Preto, SP, 15054-000, Brasil.}
%\curraddr{}
%\email{messaoud@ibilce.unesp.br}
%\thanks{The second author was partially supported by CNPq and by grant
%{\#}????, S\~ao Paulo Research Foundation (FAPESP)}



\date{} %\today

\dedicatory{}

\begin{abstract} In this paper, we  study dynamical properties as shadowing and structural stability  for a class of dynamics on $\mathbb{Z}_p$ and $\mathbb{Q}_p$, where $p \geq 2$ is a prime number. In particular, we prove that if $f: \mathbb{Z}_p \to \mathbb{Z}_p$ is a
 $(p^{-k},p^{m})$ ( $0 < m \leq k$ integers ) locally scaling map then $f$ is shadowing and  structurally stable. We also study the number of  conjugacy classes of
 these maps and we  consider the
 above properties  for $1$-Lipschitz maps of $\mathbb{Z}_p$ and for  extensions of  the shift map, contractions and dilatations on $\mathbb{Q}_p.$

\end{abstract}

\maketitle

%%%%%%%%%%%%%%%%%%%%%%%%%%%%%%%%%%%%%%%%%%%%%%%%%%%%%%%%%%%%%%%%%%%%%%%%%%%%%

\section{Introduction}

Let us consider a topological dynamical system $(X, f)$ where $(X,d)$ is a metric space  and  $f : X \to X$ is a continuous map.
A sequence $(x_n)_{n \in \mathbb{N}}$ is a {\em $\delta$-pseudotrajectory} (or a {\em $\delta$-pseudo-orbit}) of $f$
($\delta > 0$) if
$
d(f(x_n),x_{n+1}) \leq \delta $  for all nonnegative integer $ n \in \mathbb{N}.$

The map $f$ has the {\em shadowing property}
if for every $\eps > 0$ there exists $\delta > 0$ such that every
$\delta$-pseudotrajectory $(x_n)_{n \in \mathbb{N}}$ of $f$ is {\em $\eps$-shadowed}
by a real trajectory of $f$, that is, there exists $x \in X$ such that
$
d(x_n,f^n(x)) \leq \eps \text{ for all } n \in \mathbb{N}.
$

If $f: X \to X$ is an homeomorphism, the last definition is also called {\em positive shadowing}. In this case, we say that $f$ is shadowing if the last property holds  changing $\mathbb{N}$ by $\mathbb{Z}$.

The notion of shadowing, which comes from the works of
Sinai \cite{JSin72} and Bowen~\cite{RBow75},
plays a fundamental role in several branches of the area of dynamical systems.
See for instance \cite{KPal00,SPil99}.

\smallskip
Another fundamental notion in the area of Dynamical Systems is the structural stability, which was introduced by Andronov and Pontrjagin \cite{AC}.
The dynamical system $(X, f)$ is called structurally stable if there exits $\epsilon >0$ such that for all continuous map $g: X \to X$ with  $d'(f,g) <\epsilon$ ($d'$ is a distance in the space of continuous maps on $X$),
the dynamical system $(X, f)$ is topologically conjugate to the system $(X, g)$, that means that there exists an homeomorphism $h: X \to X$ such that
$f \circ h= h \circ g$.

There are other variations of the notion of  structural stability.
For instance,  a continuous map $ f : X \to X$ is said to be {\em topologically stable} if, for every $\varepsilon > 0$, there is a
$\delta > 0$ such that if $g : X \to X$ is any continuous map with $d'( f (x), g ( x ) ) < \delta$
for all $x \in X$ , then there is a continuous map (not necessarily surjective) $h : X \to X$ with
$h \circ g = f \circ h$ and $d'(h(x), x ) < \varepsilon$ for all $x \in X.$
Moreover, if $h$ is surjective, then we say that $f$ is {\em strongly topologically stable}.

It is known (\cite{W}) that if $X$ is a metric compact space and $f: X \to X$ is a shadowing expansive homeomorphism, then $f$ is topologically stable.
Moreover if $X$ is a compact topological manifold, then
the shadowing plus expansivity properties of $f$ imply
that $f$ is strongly topologically stable.
On the other hand
If $f :X \to  X $ is a topologically stable
homeomorphism of a compact manifold with topological dimension $\geq 1$,  then $f$ is shadowing (see \cite{W} for $dim (X) \geq 2$ and \cite{Mo} for $dim (X) \geq 1 $).

The last result is also true when $X$ is a Cantor compact set and was proved in \cite{K}.

Let us mention that the shadowing and stability properties were also studied in the case where $X$ is an infinite dimensional Banach space and $f: X \to X$ is a linear continuous map. See for instance \cite{BCDMP}, \cite{BM}, \cite{CGP} and \cite{Ma}.
In particular, it was proved in \cite{BM} that if $f$ is shadowing and expansive, then $f$ is hyperbolic and hence structurally stable. The converse was also investigated. In particular, the authors showed in  \cite{BM} that if $f$ is structurally stable and positively expansive, then $f$ is hyperbolic and hence it is shadowing.

In this work, we study the shadowing and  stability properties for a class of $p$-adic dynamical systems.
Let us mention that the study of dynamical systems on the ring
of p-adic integers $\Z_p$ or in its field of fractions $\mathbb{Q}_p$ is  developing (see for instance  \cite{An1}, \cite{BS}, \cite{CP}, \cite{DP}, \cite{FLY}). This study  have applications in physics, cognitive science and cryptography (see \cite{An2}, \cite{KN}).
To be precise,
let $p \geq 2$ be a prime integer and $\Z_p$ the  group of p-adic integers with its natural $p$-adic norm $\|  \; \|_p$.
We say that a map $f: \mathbb{Z}_p \to \mathbb{Z}_p$ is $(p^{-k}, p^m)$ locally scaling ( $1 \leq m \leq k$ integers) if
for all $x, y \in \mathbb{Z}_p$, whenever $\| x- y \|_p \leq p^{-k}$, we have that $\| f(x)- f(y) \|_p = p^{m} \| x- y \|_p$. This class of maps was defined in \cite{KLPS}.
An example of $(p^{-1}, p)$ locally scaling function is the shift map $S: \mathbb{Z}_p \to \mathbb{Z}_p$ defined by
$S (\sum_{i=0}^{+\infty}a_i p^{i})= \sum_{i=0}^{+\infty}a_{i+1} p^{i}.$
It is known by \cite{KLPS11}, that any $(p^{-k}, p^{k})$ locally scaling map ($k\geq 1$) is topologically conjugate to the map $S^k = S \circ  \ldots \circ S$.

In this paper, we prove that
if $f: \mathbb{Z}_p \to \mathbb{Z}_p$  is a $(p^{-k},p^{m})$ locally scaling function ($1 \leq m \leq k$), then $f$ is  shadowing and Lipschitz structurally stable (if $g$ is  closed to $f$ and $g-f$ is Lipschitz, then $g$ is topologically conjugate to $f$). As a consequence we prove that a function $f: \mathbb{Z}_p \to \mathbb{Z}_p$ is
 a $(p^{-k}, p^{k})$ locally scaling map if and only if it is topologically conjugate to $S^k$ by an isometry conjugation.
We also consider the conjugacy classes of $(p^{-k}, p^{m})$  locally  scaling maps and we prove that this set has a finite cardinality $C(m,k)$ satisfying
 $C(m, k) \geq k-m+1 $.
 Moreover $C(m,k)= k-m+1$ if and only if $m=k$. We also
prove that all $1$-Lipschitz maps on $\mathbb{Z}_p$   are shadowing and affine contractions on $\mathbb{Z}_p$   are Lipschitz structurally stable.
Let us mention that for $1$-Lipschitz  maps $f :\mathbb{Z}_p \to \mathbb{Z}_p$, the dynamical systems
$(\mathbb{Z}_p , f)$ are extensively studied (see for instance \cite{An1}, \cite{An2}, \cite{CP} and \cite{FLY}).
We also study the shadowing and stability properties for dynamics on $\mathbb{Q}_p$. In particular, we prove that for any $a, b \in \mathbb{Q}_p $, $a \neq 0$, the homeomorphism $f_{a, b}: \mathbb{Q}_p \to \mathbb{Q}_p$ defined by $f_{a,b}(z)= az+b $ for all $z \in \mathbb{Q}_p$ is shadowing. Moreover if  $\|a\|_p \neq 1$, then $f_{a,b}$ is  Lipschitz structurally stable. Moreover the conjugation is  an isometry and
$f_{a, b}$ is  topologically  conjugated to $f_{1/ p^k, 0}$ where $\|a\|_p= p^k,\; k \in \mathbb{Z} \setminus \{0\}$.
Notice that  the map $f_{1/ p^k, 0},\; k \geq 1,$ corresponds to the $k$-th iteration shift map on $\mathbb{Q}_p$.
We also prove that if $f:  \mathbb{Q}_p \to \mathbb{Q}_p$
is a dilatation or a contraction homeomorphism,
then
$f$ is shadowing and Lipschitz structurally stable.
  The last result is more general  than the previous one.  However, in the  proof, we will use the fixed point theorem to find the real orbit that shadow the $\delta$-pseudo-orbit and also to prove that the conjugation exists and is bijective. We don't know  exhibit explicitly the conjugation and we can't prove that  it is an isometry as done for  $f_{a, b}.$

The paper is organized as follows. In section 2, we consider some preliminaries and notations. In the third section, we prove the main results concerning shadowing, stability and the number of conjugacy classes of $(p^{-k}, p^m)$ locally scaling maps. In the fourth section, we consider $1$-Lipschitz maps and affine contractions on
$\Z_p$. The last section is devoted to the study of the shadowing and structural stability properties for $f_{a, b}$ on $\mathbb{Q}_p$ and also for contraction or dilatation maps on $\mathbb{Q}_p$.

\section{Preliminaries and notations}

Let $\mathbb{N}$ be the set of nonnegative integers.
The additive group of p-adic integers can be seen as the set $\mathbb{Z}_p$ of sequences $(x_i)_{i \in \mathbb{N}}$
where $x_i$ is an integer $mod \; p^{i+1}$.
Each element of  $\mathbb{Z}_p$  can be written formally as a serie
$$\sum_{i=0}^{+\infty} a_i p^i \mbox { where  } a_i \in A= \{0,\ldots, p-1 \}.$$

Let $x= \sum_{i=0}^{+\infty} a_i p^i$ and $y= \sum_{i=0}^{+\infty} b_i p^i $ be  two elements of $\mathbb{Z}_p$. We turn  $\mathbb{Z}_p$ a compact metric space by  considering the distance
$$d(x, y)= \left\{\begin{array}{cl}
0 & \mbox { if } a_i= b_i \mbox { for all } i \geq 0, \\
p^{-j} & \mbox { otherwise, where } j= min\{i \geq 0,\; a_i \ne b_i\} .
\end{array} \right. $$
We denote
$$\| x- y \|_p= d(x, y)$$  and we have $$\| x- y \|_p \leq \max \{\| x \|_p, \|  y \|_p\}.$$
We can prove that $\mathbb{Z}_p$ is the completion of $\mathbb{Z}$ under the metric $d$.

We can also add $x$ to $y$,  by adding coordinate-wise and if any of the sums is $p$ or more, we obtain the sum modulo $ p$ and we take a carry of 1 to the next sum. We can also multiply
elements of $ \Z_p$. With this, we turn $(\Z_p, +, .)$ as a ring.

The set $\mathbb{Z}_p$ can also be seen as  the projective limit of $(\mathbb{Z}/p^n \mathbb{Z}, \pi_n),\; n \geq 1 $, where for all positive integer $n$,
 $\mathbb{Z}/p^n \mathbb{Z}$ is endowed with the discrete topology, $\prod_{n=1}^{+\infty}\mathbb{Z}/p^n \mathbb{Z}$ equipped with the product topology and $ \pi_n :  \mathbb{Z}/p^{n+1} \mathbb{Z} \to  \mathbb{Z}/p^n \mathbb{Z},\;  n \geq 1,$ is the canonical homomorphism.   This limit can be seen as
  $$\mathbb{Z}_p= \{ (x_n)_{n \geq 1} \in \prod_{n=1}^{+\infty}\mathbb{Z}/p^n \mathbb{Z},\; \pi_n (x_{n+1})= x_n,\; \forall n \geq 1\}.$$
  With this definition $\mathbb{Z}_p$ is a compact set of
$\prod_{n=1}^{+\infty}\mathbb{Z}/p^n \mathbb{Z}$ (see \cite{Ro}).
 $\mathbb{Z}_p$ is a Cantor set  homeomorphic to $\Pi_{n=1}^{+\infty}A$  endowed with the  product topology of the discrete topologies on $A$.
Through this homeomorphism, any point $x= (x_n)_{n \geq 1}$ can be represented as $x= \sum_{i=0}^{+\infty} a_i p^i \mbox { where  } a_i \in A= \{0,\ldots, p-1 \}$
and $x_n = \sum_{i=0}^{n-1} a_i p^i.$

For any integer $n \geq 1, \; \mathbb{Z}/p^n \mathbb{Z}$
can be seen as a subset of $\mathbb{Z}_p$ and any element
$x $ in $\mathbb{Z}/p^n \mathbb{Z}$ can be represented  as $x=  \sum_{i=0}^{n-1} a_i p^i$ where  $a_i \in A$ for all $0 \leq i \leq n-1$.

\vspace{0.5em}

The field  of fractions of $\mathbb{Z}_p$ is denoted by $\mathbb{Q}_p$.
We have
$$ \mathbb{Q}_p= \mathbb{Z}_p [1/p]= \bigcup_{i \geq 0}p^{-i}\mathbb{Z}_p.$$
Any element  $x$ of $\mathbb{Q}_p \setminus \{0\}$ can be written as a
$$x= \sum_{i=l}^{+\infty} a_i p^i \mbox { where  } a_i \in A= \{0,\ldots, p-1 \},\; l \in \mathbb{Z} \mbox { and } a_l \neq 0.$$
We define on $ \mathbb{Q}_p$  the absolute value $\| \; \|_p$ by
$\|x\|_p= p^{-l}$ and the $p$-adic metric $d$ by
$d(x, y)= \|x - y\|_p$ for all $x, y$ in $ \mathbb{Q}_p$.

We can prove that $\mathbb{Q}_p$ is the completion of $\mathbb{Q}$ under the $p$-adic metric. Moreover  $\mathbb{Q}_p$  is not compact but locally compact set and $\mathbb{Z}_p$ is the closed unit ball of $ \mathbb{Q}_p$.

\vspace{0.5 em}
 Recall that a continuous map $f : X \to X $, where $(X, d)$ is a metric space,
is said to be {\em expansive} if there exists a constant $c > 0$ such that
for every pair $x,y$ of distinct points in $X$, there exists an integer $n \geq 0$
with $$d( f^n(x), f^n(y)) > c.$$
For example,  the shift map $S: \Z_p \to \Z_p$ is $p^{-1}$ expansive and any $(p^{-k}, p^m)$ locally scaling map ( $1 \leq m \leq k$ integers) in $\Z_p$ is $p^{-k}$ expansive.

\vspace{0.5 em}

Let $ d \geq 2$  be an integer and $\phi: A^{d} \to A$ be a function.
We say that $\phi$ is bijective  on the last variable if for all $a_1,\ldots, a_{d-1} \in A$, the function
$$\phi_{a_1, \ldots, a_{d-1}} : A \to A \mbox { defined by }
\phi_{a_1, \ldots, a_{d-1}}(x)= \phi (a_1, \ldots, a_{d-1}, x)$$
 is bijective.

\section{Shadowing and stability}

\begin{theorem}
\label{sha}

Let $f: \mathbb{Z}_p \to \mathbb{Z}_p$ be a $(p^{-k},p^{m})$ locally scaling function, where $1 \leq m \leq k $ are integers, then $f$ is  shadowing.
\end{theorem}

For the proof, we need the following Lemma.

\begin{lemma}
\label{tria}
 Let $0 < l < k$ be integers. A map $f: \mathbb{Z}_p \to \mathbb{Z}_p$ is $(p^{-k},p^{k-l})$ locally scaling if and only if there exist functions $f_i:A^{k}\longrightarrow A$ $(i=0,1,\ldots,l-1)$  and  maps $f_i:A^{k-l+i+1}\longrightarrow A$ $(\ i\geq l)$ bijective on the last variable, such that
 for all  $x=\sum_{i=0}^{\infty}x_ip^{i}\in\mathbb{Z}_p$,  we have \begin{equation}
 \label{eq}
 f(x)=\sum_{i=0}^{l-1}f_i(x_0,x_1,\ldots,x_{k-1})p^{i}+\sum_{i=l}^{\infty}f_{i}(x_0,x_1,\ldots,x_{k-l+i})p^{i}
 \end{equation}
\end{lemma}
{\bf Proof:}
Let $i \geq 0$ be an integer and
$$r_i=  \left\{\begin{array}{ll}
k-1  & \mbox  { if }  0 \leq i \leq l-1, \\
 k-l+i & \mbox { otherwise}.
\end{array}\right.$$

For all $i \geq 0$ and $x=\sum_{i=0}^{\infty}x_ip^{i}\in\mathbb{Z}_p$ , consider $$f_{i}(x_0, x_1,\ldots,x_{r_{i}})$$ as the $p^i$ coefficient of the $p$-adic representation   of $f( \sum_{i=0}^{r_{i}}x_i p^{i})$.

Since $f$ is a $(p^{-k},p^{k-l})$ locally scaling map, the functions $f_i$ are well defined. Moreover, we have that (\ref{eq}) holds and for all $i \geq l$, the functions $f_i$ are bijective on the last variable.

The converse is easy to prove and it is left to the reader.
\hfill $\Box$

\begin{remark}
\label{rpk}
It is easy to see that $f$ is a $(p^{-k}, p^k)$ locally scaling map ($k\geq 1$) if and only if there exist functions $f_i:A^{k+i+1}\longrightarrow A$ $(i\geq 0)$ bijective on the last variable, such that
$$
f(x)=\sum_{i=0}^{\infty}f_{i}(x_0,x_1,\ldots,x_{k+i})p^{i}.
 $$
  \end{remark}

As an immediate consequence of the last Lemma, we obtain the following useful result.

\begin{lemma}
\label{potencia}
 Let $0 < l < k$ be integers and $f: \mathbb{Z}_p \to \mathbb{Z}_p$ be a $(p^{-k},p^{k-l})$ locally scaling function. Then the $n$-th iterate of $x$ by $f$ is
 \begin{eqnarray}
  \label{rf}
  f^n(x)=\sum_{i=0}^{l-1}g_i^{(n)}(x_0,x_1,\ldots,x_{(n-1) (k-l)+ k-1})p^{i}+\sum_{i=l}^{\infty}g_{i}^{(n)}
  (x_0,x_1,\ldots,x_{n(k-l)+i})p^{i}
  \end{eqnarray}
  where $g_i^{(n)}:A^{n(k-l)+i+1}\longrightarrow A$ are functions bijective on the last variable, for all $i\geq l$.
\end{lemma}
%=\sum_{i=0}^{l-1}f_i(x_0,x_1,...,x_{k-1})p^{i}+\sum_{i=l}^{\infty}f_{i}(x_0,x_1,...,x_{k-l+i})p^{i}
%=\sum_{i=0}^{l-1}f_i(x_0,x_1,...,x_{k-1})p^{i}+\sum_{i=l}^{\infty}f_{i}(x_0,x_1,...,x_{k-l+i})p^{i}
\vspace{1em}
{\bf Proof:}
 We will construct the functions $g^{(n)}$ by induction on $n$.
For $n=1$, using Lemma \ref{tria}, we take
$$g_i^{(1)}=f_i, \mbox { for all } i \geq 0.$$
 Assume that for $n >1,\; g^{(n)}$ is defined and put
$$\overline{x}_{i,n}= \left\{\begin{array}{ll}
(x_0,x_1,\ldots,x_{(n-1) (k-l)+ k-1}) & \textrm{for all } 0 \leq i \leq l-1, \\
(x_0,x_1,\ldots,x_{n(k-l)+i}) & \mbox {for all } i\geq l.
\end{array}\right.$$

Define
$$g_i^{(n+1)} (x_0, \ldots, x_{n (k - l)+ k-1})=  f_i[g_0^{(n)} (\overline{x}_{0,n}),\ldots ,g_{k-1}^{(n)}(\overline{x}_{k-1, n})] \mbox {  for all } \ i=0,1, \ldots,l-1$$  and
$$g_i^{(n+1)} (x_0, \ldots ,x_{(n+1)(k-l) +i} )=f_{i}[g_0^{(n)} (\overline{x}_{0,n})\ldots,g_{k-l+i}^{(n)}(\overline{x}_{k-l+i, n})] \mbox {  for all } i \geq l.$$
 Since $f_i$ and $g_{k-l+i}^{(n)}$ are bijective functions on the last variable,  we conclude that the functions $g_i^{(n+1)}$ are also bijective  on the last variable, for all $i\geq l$. Moreover, we obtain (\ref{rf}) using Lemma \ref{tria}.
\hfill $\Box$

\begin{remark}
\label{kjk}
It is easy to see that if  $f$ is a $(p^{-k},p^{k})$ locally scaling function ($k\geq 1$),  then there exist functions $g_i^{(n)},\ i\geq 0$  from $A^{nk+i+1}$ to $A$ bijective  on the last variable, where
$$f^n(x)=\sum_{i=0}^{\infty}g_{i}^{(n)}(x_0,x_1,\ldots,x_{nk+i})p^{i}.$$
\end{remark}

\vspace{1em}

{\bf Proof of Theorem \ref{sha}}:
Assume that $m=k-l$ where $1 \leq l <k$.
Let $s$ be a nonnegative integer,
$$\varepsilon=p^{-k-s} \textrm { and } \delta=p^{-l-s}.$$
Let $(x_n)_{n \geq 0}$ a $\delta$-pseudotrajectory where $x_n=\sum_{i=0}^{\infty}x_i^{(n)}p^i$. We are going to find $y\in\mathbb{Z}_p$ satisfying
\begin{eqnarray}
\label{cond}
\|x_n - f^{n}(y)\|_p \leq p^{-k-s}, \textrm{ for all } n \geq 0.
\end{eqnarray}

Since (\ref{cond}) is valid for $n=0$, then
 \begin{eqnarray}
\label{xs}
y_i=x_i^{(0)},\ i=0,\ldots,k+s-1.
\end{eqnarray}

 Using Lemma \ref{tria}, we have
$$f(y)=\sum_{i=0}^{l-1}f_i(y_0,y_1,\ldots,y_{k-1})p^{i}+\sum_{i=l}^{\infty}f_{i}(y_0,y_1,\ldots,y_{k-l+i})p^{i}.$$

Since (\ref{cond}) is valid for $n=1$, it follows that
\begin{equation}
\label{r1}
f_i(y_0,y_1,\ldots,y_{k-1})=x_i^{(1)}, \textrm{ for } i=0,\ldots,l-1,
\end{equation}
and
\begin{equation}
\label{r2}
f_i(y_0,y_1,\ldots,y_{k-l+i})=x_i^{(1)}, \textrm{ for } i=l,\ldots,k+s-1.
\end{equation}

Since  $$\|x_1 - f(x_0)\|_p \leq p^{-l-s},$$ by (\ref{xs}) we deduce that the equations (\ref{r1}) and (\ref{r2}) are satisfied. Since $y_0,\ldots,y_{k+s-1}$ have been  determined and $f_{l+s}$ is  bijective on the last variable, we can find  $y_{k+s} \in A$ such that $$f_{l+s}(y_0,y_1,\ldots,y_{k+s-1}, y_{k+s})=x_{l+s}^{(1)}.$$
Continuing by the same way and using recursively the fact that
$$f_{l+s+1}, \ldots, f_{k+s-1}$$
 are bijective on the last variable, we find recursively
$$y_{k+s+1},\ldots,y_{2k+s-l-1} \in A \mbox { satisfying the
relation } (\ref{r2}).$$

%$$f(y)=\sum_{i=0}^{l-1}f_i(x_0^{(0)},x_1^{(0)},...,x_{k-1}^{(0)})p^{i}+\sum_{i=l}^{ %\infty}f_{i}(x_0^{(0)},x_1^{(0)},...,x_{k-1}^{(0)},y_k,...,,y_{k+i})p^{i}.$$
Let $n>1$ and suppose, by induction, that we have found
 $$y_0,y_1,\ldots,y_{n k-(n-1)l+s-1}.$$

By (\ref{cond}) and Lemma \ref{potencia}, we have
$$f^n(y)=\sum_{i=0}^{k-1}x_i^{(n)}p^{i}+\sum_{i=k}^{\infty}g_{i}^{(n)}(y_0,y_1,\ldots,y_{n(k-l)+i})p^{i}.$$

Since $$\|x_{n+1}- f(x_n)\|_p \leq p^{-l-s},$$
then $$x_i^{(n+1)}=f_i(x_0^{(n)},x_1^{(n)},\ldots,x_{k-1}^{(n)}),\ i=0,1,\ldots,l-1.$$

On the other hand

\smallskip

\noindent $f^{n+1}(y)= f(f^n(y))$

\hfill $\displaystyle=\sum_{i=0}^{l-1}f_i(x_0^{(n)},x_1^{(n)},\ldots,x_{k-1}^{(n)})p^{i}+\sum_{i=l}^{\infty}f_{i}(x_0^{(n)},x_1^{(n)},\ldots,x_{k-1}^{(n)},g_k^{n}(\overline{y}_{k, n}),\ldots,g_{k-l+i}^{n}(\overline{y}_{k-l+i, n}))p^{i}$

\smallskip

where
$$\overline{y}_{j, n}= (y_0,y_1,\ldots,y_{n(k-l)+j}), \textrm{ for all }  j=k,\ldots, k-l+i.$$

Since (\ref{cond}) holds for $n+1$, we have
$$f_{i}(x_0^{(n)},x_1^{(n)},\ldots,x_{k-1}^{(n)},g_k^{n}(\overline{y}_{k, n}),\ldots,g_{k-l+i}^{n}(\overline{y}_{k-l+i, n}))=x_i^{(n+1)},i=l,\ldots,k+s-1.$$

Since the functions $f_i$ and $g_i^n\ (i\geq l)$ are bijective on the last variable, we can find recursively
$$y_{nk-(n-1)l+s},\ldots, y_{(n+1)k-nl+s-1}.$$

The last proof can be done by the same manner when $m=k$, using  Remarks \ref{rpk} and \ref{kjk}. The unique difference in this case is that
 $$ \mbox { if } \varepsilon=p^{-k-s} \mbox {  then } \delta=p^{-k-s}.$$
\hfill$\Box$

\begin{theorem}
\label{strucut}

Let  $1 \leq m \leq k $  be integers. If $f: \mathbb{Z}_p \to \mathbb{Z}_p$ is a $(p^{-k},p^{m})$ locally scaling map, then $f$ is Lipschitz structurally stable.
\end{theorem}

\vspace{0.5em}

\begin{remark}
If $f: \mathbb{Z}_p \to \mathbb{Z}_p$ is $(p^{-k},p^{m})$ locally scaling map $(1 \leq m \leq k $ integers), then $f$ can not be structurally stable.

Indeed, for all $j \geq 0$, define $$\phi_j (x)= \sum_{i=0}^{j} f(x)_i p^i$$ where $f(x)_i$ is the $p^i$ coefficient of $f(x)$ in its $p$-adic representation.
Assume that $h$ is a conjugation between $f$ and $\phi_j$ ($j$ large).  Thus for all $n >j$,
$$h(f^n (x))= \sum_{i=0}^{j} f(h(x))_i p^i \mbox { for all } x \in \Z_p.$$
Hence $h$ is not injective.
\end{remark}

\vspace{0.5em}

\noindent {\bf Notation}:
Let $f :\mathbb{Z}_p \to \mathbb{Z}_p$ be a shadowing map and  $\varepsilon >0$. We denote $\delta_{f}(\varepsilon)$ as the supremum of $\delta>0$ such that any $\delta$-pseudo-orbit of $f$ can be $\varepsilon$ shadowed by a real orbit.

\begin{remark}
\label{rfg}
By the proof of Theorem 1, we deduce that if $f$ is $(p^{-k},p^{m})$ locally scaling map where $  1 \leq m \leq k $ are integers,  then
$$p^{m-k} \leq \delta_{f}(p^{-k}) \mbox { if } m <k$$
 and
 $$p^{-k} \leq \delta_{f}(p^{-k}) \mbox { if } m=k.$$
 \end{remark}

\begin{lemma}
\label{str}

Let $ k \geq 1$ be an integer and  $f, g: \mathbb{Z}_p \to \mathbb{Z}_p$  be  $p^{-k}$  expansive and shadowing maps. If
 $$\|f - g\|_{\infty} < \delta= min (\delta_{f}(p^{-k}),\delta_{g}(p^{-k})),$$ then $f$ and $g$ are topologically conjugate.
\end{lemma}

\begin{remark}
The proof of Lemma \ref{str} is similar to the proof of the result
that states that a shadowing expansive map on a $d \geq 1$ dimensional compact manifold is topologically stable (see \cite{SPil99} and \cite{W}).
\end{remark}
{\bf Proof of Lemma \ref{str}:}
Let $x\in\mathbb{Z}_p$, then we have  $$\| g^{n+1}(x)-f(g^n(x))\|_p=\| (g-f)(g^n(x))\|_p < \delta.$$
Since $f$ is shadowing, there exists  $y= h(x) \in\mathbb{Z}_p$ such that, for all $n\in\mathbb{N}$,
\begin{eqnarray}
\label{f11}
\| g^n(x)-f^n(h(x))\|_p \leq p^{-k}.
\end{eqnarray}

Since $f$ is $p^{-k}$ expansive, we deduce that
$h(x)$ is unique.

\smallskip

\textbf{Claim 1:} $f\circ h=h\circ g$.

\smallskip

We have that
 $$\| g^{n+1}(x) - f^n(f\circ h(x)) \|_p \leq p^{-k} \mbox {  and }
\| g^{n}(g(x)) - f^{n}( h(g(x))) \|_p \leq p^{-k}, \textrm{ for all }   n \geq 0.$$
 Thus

\begin{center}
$\| f^n(f\circ h(x))-f^n(h\circ g(x))\|_p \leq p^{-k}$, for all $n\in\mathbb{N}$.
\end{center}
Since $f$ is $p^{-k}$ expansive, it follows that
 $$f\circ h(x)= h\circ g(x).$$

\smallskip

\textbf{Claim 2:} $h$ is bijective.

Given $x \in \mathbb{Z}_p$ and changing the roles of $f$ and $g$, we deduce as done before that  there exists a unique $ k(x) \in\mathbb{Z}_p$ such that for all $n\in\mathbb{N}$
\begin{eqnarray}
\label{f2}
\| f^n(x)-g^n(k(x))\|_p \leq p^{-k}.
\end{eqnarray}

By(\ref{f11}) and (\ref{f2}), we deduce that
for all integer $n \geq 0$ and $x \in \mathbb{Z}_p$,
$$ \|g^{n} (x)- g^n (k \circ h(x))\| \leq p^{-k},\; \|f^{n} (x)- f^n (h \circ k(x))\| \leq p^{-k}.$$
Hence $$k \circ h(x)= h \circ k(x)=x$$
and we are done.

\smallskip

\textbf{Claim 3:} $h$ is continuous.

\smallskip

Let $\gamma >0$. Since $\Z_p$ is compact, there exists an integer $N>0$ sufficiently large such that for all $z, w \in
\Z_p$,
if
$$\|f^n (z)- f^n (w) \|_p < p^{-k} \mbox { for all } n=0,\ldots, N,$$
then
$$\|z-w\|_p <\gamma.$$

Let $x, y \in \Z_p$. There exists a real number $\alpha >0$ such that if $\|x-y\|_p \leq \alpha$, then
\begin{eqnarray}
\label{dfbl1}
\|g^n(x)- g^n(y) \|_p \leq p^{-k}, \textrm{ for all } n=0,\ldots, N.
\end{eqnarray}
Hence if $\|x-y\|_p \leq \alpha$, then we deduce
by (\ref{f11}) and (\ref{dfbl1}) that
\begin{eqnarray}
\|f^n(h(x))- f^n(h(y)) \|_p \leq p^{-k}, \textrm{ for all } n=0,\ldots, N.
\end{eqnarray}
Thus $$\|h(x)- h(y) \|_p < \gamma$$ and we obtain the claim.
\hfill $\Box$

\vspace{1em}

{\bf Proof of Theorem \ref{strucut}:}
Since $f$ is $(p^{-k}, p^{k-l})$ locally scaling, then $f$ is $p^{-k}$ expansive. It is also shadowing by Theorem \ref{sha}.
Let $g= f+\psi$ where
$\psi : \mathbb{Z}_p \to \mathbb{Z}_p$  is a
$$ \delta <  p^{k-l-1}(p-1) \mbox { Lipschitz map }.$$

{\bf Claim}: $g$ is $(p^{-k}, p^{k-l})$ locally scaling.

Indeed,
let $ x, y \in \mathbb{Z}_p$ such that $\| x - y \| \leq p^{-k}$.
We have
$$
\begin{array}{rcl}
(p^{k-l}- \delta) \|x- y\|_p & \leq & \|g(x)- g(y)\|_p \\
                           & \leq & \max\{p^{k-l} \|x- y\|_p,\delta \|x- y\|_p\} \\
                           & = &  p^{k-l} \|x- y\|_p.
\end{array}
$$
We deduce that
$$\|g(x)- g(y)\|_p = p^{k-l} \|x- y\|_p $$
and we are done.

Using the claim, we deduce that $g$ is $p^{-k}$ expansive and shadowing map. By  Lemma \ref{str}, we are done.
\hfill $\Box$

Another  consequence is  the following  result.

\begin{proposition} \label{propiso}
Let $f: \mathbb{Z}_p \to \mathbb{Z}_p $ be a function and $k$ be a positive integer. Then $f$ is a
 $(p^{-k}, p^k)$ locally scaling map if and only if $f$ is topologically conjugate by an isometry conjugation to the map $S^k$ where $S$ is the shift map.
\end{proposition}

\vspace{1em}

\begin{remark}
 The fact that a $(p^{-k}, p^k)$ locally scaling map  is topologically conjugate to the map $S^k$ was first proved in \cite{KLPS11}. Here, we give another proof using the shadowing property.
\end{remark}

\vspace{0.5 em}

\begin{remark}
If
$$\|S^k -f \|_{\infty} \leq p^{-k},$$
 then by Remark \ref{rfg} and Lemma \ref{str},
we deduce that $f$ is conjugate to the map $S^k$.
\end{remark}

\vspace{0.5 em}
{\bf Proof of Proposition \ref{propiso}}: The converse is easy to prove and it is left to the reader.

Now assume that $f$ is $(p^{-k}, p^k)$ locally scaling. For all $x= \sum_{i=0}^{+\infty} x_i p^i$, put $f(x)=
\sum_{i=0}^{+\infty} f(x)_{i} \; p^i$ and
$h(x)=
\sum_{i=0}^{+\infty} h_{i}(x) p^i$ where
$$h_i (x)= f^q(x)_r \mbox { when } i= qk+r,\; q, r \in \mathbb{N} \mbox { and }0 \leq r \leq k-1.$$
That is
\begin{eqnarray}
\label{hnm}
h_0(x)h_1(x) \ldots= x_0 \ldots x_{k-1}\; f (x)_0 \ldots f (x)_{k-1} \ldots f^i (x)_0 \ldots f^i (x)_{k-1} \ldots .
\end{eqnarray}

For all nonnegative integer $n$, we have that
\begin{eqnarray}
\label{gg1}
\| f^n(x)-S^{nk}(h(x))\|_p \leq p^{-k}.
\end{eqnarray}

Observe that as done in claim 1 of the proof of Lemma \ref{str}, we have
$$S^k \circ h= h \circ f.$$

{\bf Claim:} $h$ is  an isometry.

Indeed, let $x, y \in \mathbb{Z}_p$ such that $\|x- y\|_p = p^{-m},\; m \geq 0$. Put
$$m= qk+r,\; q, r \in \mathbb{N} \mbox { and } 0 \leq r \leq k-1.$$
 Since $f$ is $(p^{-k}, p^{k})$ locally scaling, we have
$$\|f^i(x)- f^i(y)\|_p=  p^{ik-m}, \textrm{ for all } 0 \leq i \leq q.$$
Hence, for all  $0 \leq i \leq q$, $f^{i}(x)$ and $f^{i}(y)$ coincide exactly on $m-ik$ first digits. Therefore, we have that
$$\|h(x)- h(y)\|_{p} = p^{-m}= \|x- y\|_p.$$

To prove that $h$ is an homeomorphism, it suffices to prove that $h$ is surjective.

Indeed let
 $y= \sum_{i=0}^{+\infty} y_i p^i \in \mathbb{Z}_p$.
 We aim to find
 $$x= \sum_{i=0}^{+\infty} y_i p^i \in \mathbb{Z}_p \mbox { such that } h(x)= y.$$
 By (\ref{hnm}), we deduce that

 $$x_i= y_i, \textrm{ for all } 0 \leq i \leq k-1.$$

 Assume that we have constructed $$x_i, \textrm{ for all } 0 \leq i \leq nk-1 \mbox { where } n \geq 1 \mbox { is an integer }.$$

 Since $f$ is $(p^{-k}, p^k)$ locally scaling, using Remark \ref{kjk}, we deduce that
 $$y_{nk+i}= f^n(x)_i=
  g_{i}^{(n)}(x_0,x_1,\ldots,x_{nk+i}), \textrm{ for all } 0 \leq i \leq k-1$$  where the functions $g_i^{(n)}$ are functions from $A^{nk+i+1}$ to $A$ bijective  on the last variable.
 We obtain recursively
 $$x_{nk},\ldots, x_{(n+1)k-1}$$
 and we are done
 \hfill $\Box$

 \begin{definition}
Let $1 \leq m \leq k$ be integers. A set $S(m,k)$ of conjugacy classes of
 $(p^{-k}, p^{m})  $ locally scaling maps is by definition a subset of $(p^{-k}, p^m)$   locally scaling maps such that
 any $(p^{-k}, p^{m})  $ locally scaling map is topologically conjugate to some element of  $S(m,k)$  and any two different elements of   $S(m,k)$  are not topologically conjugate.

 \end{definition}

 \begin{theorem}\label{card}
Let $1 \leq m \leq k$ be integers. Any set of conjugacy classes of
 $(p^{-k}, p^{m})  $ locally scaling maps has a finite cardinality $C(m,k) \geq k-m+1$.
 Moreover $C(m,k)= k-m+1$ if and only if $m=k$.
\end{theorem}

\begin{remark}
By Proposition \ref{propiso}, we deduce that $C(k, k)= 1$ for all integer $k \geq 1$.
\end{remark}

{\bf Proof of Theorem \ref{card}:}
First observe that any two sets of conjugacy classes of
 $(p^{-k}, p^{m})  $ locally scaling maps have the same   cardinality $C(m,k)$.
 Now, let $T_i,\; i=0,\ldots, k-m$, be the  functions in $\Z_p$ defined by
$$T_0= S^m $$ where $S$ is the shift map
and for all $j=1,\ldots,k-m$ and $x= \sum_{i=0}^{+\infty} x_i p^i$ in $\Z_p$,
$$T_j\left(\sum_{i=0}^{+\infty} x_i p^i\right)=\sum_{i=0}^{+\infty} y_i p^i$$ where
$$y_0 y_1 y_2 \ldots= x_0 \ldots x_{j-1} x_{m+j}x_{m+j+1}\ldots$$
Observe that
$T_0$ is a $(p^{-m},p^m)$   locally scaling function and therefore a $(p^{-k},p^m)$ locally scaling function. On the other hand, by using Lemma \ref{tria}, we deduce that  the maps $T_j$ are  $(p^{-m-j},p^m)$ locally scaling function and hence $(p^{-k},p^m)$ locally scaling function.
Now we are going to calculate the number of fixed points of these functions. We can show that
$$S^{m}(x)=x \textrm{ if and only if } x=\sum_{t=0}^{+\infty}(x_0p^{tm}+x_1p^{tm+1}+...+x_{m-1}p^{tm+m-1})$$ and so this function has $p^{m}$ fixed points. We can also show that

$$T_j(x)=x \textrm{ if and only if } x=\sum_{i=0}^{j-1}x_ip^{i}+\sum_{t=0}^{+\infty}(x_jp^{tm+j}+x_{j+1}p^{tm+j+1}+...+x_{m+j-1}p^{tm+m+j-1})$$  and so these functions have $p^{m+j}$ fixed points.\\ Since the number of fixed points are different for these functions then they are not pairwise conjugate and then we have $k-m+1$ different conjugacy classes. Thus
 $$C(m,k) \geq k-m+1.$$

 \vspace{0.5em}

   Now assume that $m <k$ and
   consider the map $$R(x)= S^m(x) \mbox { if } x= \sum_{i=0}^{\infty}x_{i} p^{i},$$
    where $x_0 \in \{0,\ldots, p-2\}$
   and $$R(x)= T_1(x) \mbox { if } x_0= p-1.$$

It is easy to see that $R$ is a $(p^{-m-1},p^m)$ locally scaling function and hence $(p^{-k},p^m)$ locally scaling function. On the other hand $R$ has
   $$p^{m-1} (p-1)+ p^{m+1}$$
    fixed points which are different for $p^{m+j}$ for all $0 \leq j \leq k-m$.
   Hence $R$ cannot be conjugate to some $T_i,\; i \in \{0,\ldots, k-m\}$. Thus
   $$C(m,k) > k-m+1.$$

    \vspace{0.5em}

{\bf Claim}: For any integers $1 \leq m \leq k$, the number $ C(m, k)$ is finite.

Indeed, let $F(m,k)$ be the set of $(p^{-k},p^m)$ locally scaling maps in $\Z_p$.
We have $$F(m,k) \subset F(m,k+1).$$
Now let $f \in F(m,k+1) \setminus F(m,k) $.

Then by Lemma \ref{tria},
  for all  $x=\sum_{i=0}^{\infty}x_{i} p^{i}\in\mathbb{Z}_p$,  we have $$
 f(x)=\sum_{i=0}^{k-m}f_i(x_0,x_1,\ldots,x_{k})p^{i}+\sum_{i=k-m+1}^{\infty}f_{i}(x_0,x_1,\ldots,x_{m+i})p^{i}
$$
   where  $f_i,\ i=0,1,\ldots,k-m,$ are functions from $A^{k+1}$ to $A$ and $f_i,\ i\geq k+1-m$, are functions from $A^{k-m+i+2}$ to $A$ bijective on the last variable. Moreover,
   there exists an integer $0 \leq i \leq k-m-1$ such that
   $f_i$ depends also of $x_k$ or $f_{k-m}$ is not bijective in the last variable.
   Observe that there exists a finite number of $f_i,\; 0 \leq i \leq k-m-1$,
    satisfying the above property. Moreover, if $g$ is a $(p^{-k-1}, p^m)$ locally scaling map such that $\| f(z)- g(z) \|_p \leq p^{m-k}$ for all $z\in\mathbb{Z}_p$, then by Lemma \ref{str}, $g$ is topologically conjugate to $f$. Hence, there exists a finite number of functions $f$ belonging to $F(m,k+1) \setminus F(m,k) $. Therefore, the claim is proved.

   \hfill $\Box$

   \begin{remark}
   1.
If
$\|f- T_i \|_{\infty} \leq p^{-l},$ for some  $i \in \{0,\ldots k-m\}$,
 then by Remark \ref{rfg} and Lemma \ref{str},
we deduce that $f$ is topologically conjugate to the map $T_i$.

\noindent 2.   It will be interesting to characterize all  conjugacy classes of the set of
 $(p^{-k}, p^{m}) $ locally scaling maps.

 \end{remark}

\section{One Lipschitz maps on $\mathbb{Z}_p$}

\begin{proposition}
\label{lipc}
If $f: \Z_p \to \Z_p$ is a $1$-Lipschitz map, then $f$ is shadowing.
\end{proposition}

{\bf Proof:} Let $\varepsilon >0$  and $(x_n)_{n \geq 0}$ be a {\em $\varepsilon$-pseudotrajectory} of $f$.

Put $$w_n= x_{n+1}- f (x_n),\; \mbox { for all } n \geq 0.$$
Thus ,
$$ x_{n+1} - f^{n+1} (x_0) = f(x_n) - f^{n+1} (x_0) + w_n, \; \mbox { for all } n \geq 0.$$

\vspace{1em}

Hence
$$ \| x_{n+1} - f^{n+1} (x_0)\|_p  \leq \max \{  \| x_n - f^{n} (x_0)\|_p , \vert w_n \vert_p\}.$$
We deduce by induction that
 for all $n \geq 0$,
 $$
  \| x_{n+1} - f^{n+1} (x_0)\|_p  \leq \max \{  \| w_0 \|_p ,  \| w_1 \vert_p, \ldots,\vert w_n \|_p\} \leq \varepsilon .
$$
\hfill $\Box$

\begin{remark}
\label{eza}
If  $f: \mathbb{Q}_p \to \mathbb{Q}_p$ is a  $1$-Lipschitz function, then $f$ is positively shadowing.
\end{remark}

\vspace{0.5 em}

\begin{corollary}
If  $f: \Z_p \to \Z_p$ is an analytic  map (in particular a polynomial map), then $f$ is shadowing.
\end{corollary}

\vspace{0.5em}

Another consequence of Proposition \ref{lipc} is the following.

Given $A= \{0,\ldots, p-1\}$ and $A^{*}$ be the set of finite words on the alphabet $A$. Let $\sigma: A \to A^{*}$ be a map such that for any $a \in A$,   $\sigma (a)$ is not the empty word. The map $\sigma$ is called substitution.
Now let $A^{\N}$ be the set of infinite words on $A$.
We extend $\sigma$ to $A^{\N}$ by concatenation and denote it also by $\sigma$. That is
$$\sigma (a_0 a_1 \ldots)= \sigma (a_0)\sigma (a_1) \ldots, \mbox { for all } a_0 a_1 \ldots \in A^{\N}.$$

Let us mention that substitutions and their associated symbolic dynamical systems have been extensively studied (see for instance \cite{PF} and \cite{Q}).

Now, given a substitution $\sigma$ on $A$, we define the map
$\sigma: \Z_p \to Z_p$  (called substitution on $\Z_p$)
by
$$\sigma ( \sum_{i=0}^{+\infty} a_i p^i)= \sum_{i=0}^{+\infty} b_i p^i$$
where
$$ \sigma (a_0 a_1 \ldots)= b_0 b_1 \ldots.$$
Since $\sigma$ is $1$-Lipschitz on $\Z_p$, we deduce that

\begin{corollary}

Any substitution on $\Z_p$ is shadowing.

\end{corollary}

\vspace{0.5em}

\noindent {\bf Question:} Is a substitution on $\Z_p$ can be Lipschitz structurally stable ?

\vspace{0.5em}

It is known (\cite{Ma1}, \cite{Ma2}) that any continuous map $f: \Z_p \to \Z_p$ can be represented by the Mahler series
$$ f(x)= \sum_{n=0}^{+ \infty} a_{n} \left(\begin{array}{c} x \\ n \end{array}\right) $$
where $a_n \in \Z_p$ for all integer $n \geq 0$ and the binomial coefficient functions are defined by
$$\left(\begin{array}{c} x \\ n \end{array}\right)= \frac{1}{n!} x (x-1) \ldots (x-n+1) \; (n \geq 1)  \mbox { and }  \left(\begin{array}{c} x \\ 0 \end{array}\right)=1.$$
Moreover,
the function $f$ is $1$-Lipschitz if and only if (see
\cite{An1} and \cite {An2}) for all $n \geq 0$,
 $$\| a_n \|_p \leq p^{-[\frac{log \; m} {log \; p}]}.$$

If the function $f$ is equal to  $S^k$, then (see \cite{KLPS11})
$$a_n= 0, \mbox { for all }  0 \leq n <p^k  \mbox { and } a_{p^k}=1  $$
and

$$ \| a_n \|_p \leq p^{-j} \mbox { for } n > j p^{k}- j+1.$$

\vspace{1em}

\noindent {\bf Question:}
Can we characterize the shadowing property of a map $f$ by the coefficients of its Mahler series?

\vspace{1em}

\begin{remark}
\label{rgm}
Let $X \in \{\Z_p,\; \mathbb{Q}_p$\} and $f: X \to X$ be a  $1$-Lipschitz map. Then $f$ is not necessarily Lipschitz topologically  stable.
For example, the function identity $id$  is shadowing, since it is $1$-Lipschitz but not Lipschitz topologically  stable, since if we consider the function $g$ defined by
$$g(z)= z+ p^{k} z,\; k \mbox { large enough }$$ and $h:X \to X$ such that  $h \circ id= g \circ h$, we obtain $h=0$.

It will be interesting to verify if there are $1$-Lipschitz maps that are Lipschitz topologically  stable.
\end{remark}

\vspace{1em}

\noindent {\bf Question:}
 Assume that $f: \Z_p \to \Z_p$ is shadowing and $c$ expansive where $0< c \leq p^{-1}$.
Is $f$ Lipschitz structurally stable?

\begin{proposition}
Let $a, b$ be elements of $\mathbb{Z}_p$ such that $\|a\|_p <1$.  Let $f: \mathbb{Z}_p \to \mathbb{Z}_p$ be the affine map  defined by $f(z)= a z+b$ for all  $z \in \mathbb{Z}_p$, then $f$ is Lipschitz structurally stable.
\end{proposition}

Observe that $f$ is a contraction and is not surjective.

\vspace{1em}

{\bf Proof:}
First assume that $a= p$.
Since $f$ is topologically conjugate to the map $f_0$ defined by $f_0 (z)= pz$ for all $z \in \mathbb{Z}_p$ by the isometry $\phi: \mathbb{Z}_p \to \mathbb{Z}_p$ defined by $\phi(z)= z + \frac{b}{1- p} $, we can
 suppose that $b=0$.
 Let  $g: \mathbb{Z}_p \to \mathbb{Z}_p$
such that $g- f= \psi$ is a $\delta <1/p$ Lipschitz map and $\| g- f\|_{\infty} \leq \delta$.
Let us construct  $ h: \mathbb{Z}_p \to \mathbb{Z}_p$ such that
\begin{eqnarray}
\label{sxan}
 h \circ f= g \circ h.
 \end{eqnarray}

(\ref{sxan}) is equivalent to
\begin{eqnarray}
\label {azx}
h(p z)= p \; h(z)+ \psi (h(z)),\; \forall z \in \mathbb{Z}_p.
\end{eqnarray}

Define $$h(z)= z \mbox { for all } z \in E_0= \mathbb{Z}_p \setminus p \mathbb{Z}_p.$$
For all $n \geq 1$, we define recursively, using (\ref{azx}), $h$ on $ E_n= p ^n\mathbb{Z}_p \setminus p^{n+1} \mathbb{Z}_p.$
Precisely,
\begin{eqnarray}
\label {qwbn}
h(p^{n} z)= p^n z+ p^{n-1} \psi(z)+ \sum_{i=0}^{n-2} p^i \psi \circ h \; (p^{n-1-i}( z)) \mbox { for all } n \geq 1 \mbox { and } z \in E_0.
\end{eqnarray}

Notice that that all $E_n$ are disjoint and
$\mathbb{Z}_p= \{0\} \cup \bigcup_{n=0}^{+\infty} E_n$. We put $h(0)=0$, then
 $g(0)= \psi(0)=0.$

\vspace{1em}

{\bf Claim 1}: For all integer $n \geq 0$, $h(E_n) = E_n$ and $h: E_n \to E_n$ is continuous and bijective.

 \vspace{1em}
 Indeed, if $n=0$, it is clear.
 Assume that  $h(E_n) =  E_n$ for some integer $n \geq 0$
  and $h: E_n \to E_n$ is continuous and bijective.

 Let $y \in E_n$ and let
 $$r: E_n \to E_n,\; \textrm { defined by }  r(z)= y - \frac{ \psi (z)}{p}.$$
 Since $E_n$ is a closed set on $\mathbb{Z}_p$ and $r$ is a contraction map, we deduce that there exists $z'$ in $E_n$ such that
 $ p y= \psi (z')+ p z'.$
 Hence, by induction hypothesis, we deduce that
 there exists $z \in E_n$ such that
 $$p y = p h(z)+ \psi (h(z))= h(pz).$$
 Hence, we have that $E_{n+1} \subset h(E_{n+1}).$

 On the other hand, since $\psi(0)=h(0)=0$ and $\psi$ is $\delta$ Lipschitz, we deduce by (\ref{azx}) that
 \begin{eqnarray}
 \label{lsc}
 \| \psi (h(z)) \|_p \leq \delta \|h(z)\|_p = \delta /p^{n} < 1/p^{n+1}=\| p h(z)\|_p , \; \textrm  { for all } z \in E_n.
 \end{eqnarray}
Thus, using (\ref{azx}) and (\ref{lsc}), we deduce that if $z \in E_n$, then
 $$\|h(p z)\|_p =  max (\| p h(z)\|_p, \| \psi (h(z))\|_p)  =1/p^{n+1}.$$
 Hence $h(pz) \in E_{n+1}$ and we deduce that
 $E_{n+1} = h(E_{n+1}).$

 Now let  us prove that $h / E_{n+1}: E_{n+1} \to E_{n+1}$ is injective.

 Indeed, by (\ref{azx}), we deduce that
 $$\| h(pz)- h(pz')\| \geq (\frac{1}{p}- \delta) \| h(z)- h(z')\|_p,\; \textrm{ for all } z, z' \in E_{n}.$$

 Since $h: E_{n} \to E_{n}$ is injective, then $h / E_{n+1}$ is too.

 \vspace{1em}

 {\bf Claim 2}: $h$ is continuous on $0$.

 \vspace{1em}

Indeed, using (\ref{qwbn}), we deduce
that for all $z$ in $ E_0$ and $n $ positive integer
$$ \| h(p^{n} z)\|_p \leq max\{ p^{-n}, p^{-i} \|\psi \circ h \; (p^{n-1-i}( z)) \|_p,\; i=0,\ldots, n-2\} .$$
Since $\|\psi \circ h \; (p^{n-1-i}( z)) \|_p \leq \delta  \| h  (p^{n-1-i}( z)) \|_p < p^{-n+i}$,
we deduce that $\| h(p^{n} z)\|_p \leq  p^{-n}$ and we are done.

\vspace{1em}

Finally, if $\|a\|_p = p^{-k},\; k \geq 1$, the proof is the same than below.
 Since $h(a z)= a h(z)+ \psi (h(z))$ for any $z \in \mathbb{Z}_p$, we choose $E_n = a^{n} \mathbb{Z}_p \setminus
 a^{n+1 }\mathbb{Z}_p$ for all $n \geq 0$
 and we define $h(z)= z$ for all $z\in E_0$.

 \hfill $\Box$

 \begin{remark}
 The conjugation $h$ defined in the last proposition is not unique.
 \end{remark}

 \vspace{1em}

 {\bf Question:} Is any contraction on $\mathbb{Z}_p$ is Lipschitz structurally stable?
 We will see in the next section  that this is true for contractions on $\mathbb{Q}_p.$

\vspace{1em}

\begin{remark}
There exists a class of non shadowing maps on $ \mathbb{Z}_p$.

Indeed,  let $B$ be a finite alphabet and $B^{\N}$ be the set of infinite words endowed with the topology product of discrete topologies on $B$.
Now consider $X$ as a closed subset of $B^{\N}$. Assume that $X$ is a subshift of $B^{\N}$. That means $S(X)$ is contained in $X$ where $S$ is the shift map.
It is classical (see \cite{Aok89}) that $S: X \to X$ is shadowing if
and only if $(X,S)$ is of finite type.

Now, consider $(X, S)$ a subshift which is not of finite type. Then it is not shadowing.

Since $X$ is Cantor set, there exists an homeomorphism $\phi$ from $X$ to $\Z_p$. Define the map $f$ on $\Z_p$ by
$f \circ \phi= \phi \circ S.$ Since $S$ is not shadowing, then $f$ is not shadowing.
\end{remark}

\section{Shadowing and stability on $\mathbb{Q}_p$}

Let $a \in \mathbb{Q}_p$ and
$g_a: \Z_p \to \Z_p$  be the map defined by
$$g_a(z)= az \mbox{ mod } (\Z_p) \mbox { for all }  z \in \Z_p$$
where  for all $x = \sum_{i=k}^{\infty} a_i p^{i} \in \mathbb{Q}_p,\; k \leq 0$, we have
$x \mbox{ mod } (\Z_p)= \sum_{i=0}^{\infty} a_i p^{i}$.

Observe that if $a=1/p$, then $g_a$ is equal to the shift map $S$.

In \cite{KLPS11}, the authors showed that if $\| a \|_p= p^k,\; k>0$, then $f_a$ is $(p^{-k}, p^k)$  locally scaling, and hence conjugate  to $S^k.$

 \vspace{1em}

 Now let us consider   $a, b \in \mathbb{Q}_p  $ such that $b \neq 0$ and
 $f_{a, b}: \mathbb{Q}_p \to \mathbb{Q}_p$ defined by:
 $$f_{a, b}(z)= a z+ b.$$

 Observe that $f_{a, b}$ is an homeomorphism and $f_{1/p, 0}$ corresponds to the shift map on $\mathbb{Q}_p$.

  \vspace{1em}

\begin{theorem}
\label{thb}
 Let  $a, b $ in  $\mathbb{Q}_p $ such that $a \neq 0$. Then $f_{a,b} : \mathbb{Q}_p \to \mathbb{Q}_p$ is shadowing. Moreover, if $\|a\|_p \neq 1$, then $f_{a,b}$ is Lipschitz structurally stable.  Furthermore, the conjugation is an isometry and
$f_{a, b}$ is  topologically  conjugate to $f_{ 1/ p^k, 0}$, where $\|a \|_p= p^k,\; k \in \mathbb{Z}\setminus \{0\}$.
 \end{theorem}

 \begin{lemma}
 \label{swvb}
Let $m \in \mathbb{Z}$ be an integer and  $\phi:\mathbb{Q}_p \to \mathbb{Q}_p $ be a map such that
$$\|\phi(x)- \phi(y)\|_p= p^m \|x- y\|_p \mbox { for all } x, y \mbox { in } \mathbb{Q}_p.$$
 Then there exist functions $\phi_i,\ i\in \mathbb{Z}$,  from $\prod_{n= -\infty}^{i+m} A_n$ to $A$  bijective on the last variable, where $A_n= A$ for all integer $n$ and $\phi_i (\ldots 0 \ldots 0)= 0$ for all $i \in \mathbb{Z}$, such that
for all $x=\sum_{i=-\infty}^{+\infty}x_ip^{i}\in\mathbb{Q}_p$, we have
\begin{eqnarray}
\label{fgb}
\phi(x)=\sum_{i=-\infty}^{+\infty}\phi_{i}(\ldots x_{i+m-1} x_{i+m})p^{i}.
 \end{eqnarray}

 \end {lemma}

 { \bf Proof}:
 For all $i $ in $\mathbb{Z}$ and $x=\sum_{j=-\infty}^{\infty}x_j p^{j}\in\mathbb{Q}_p$ , consider $\phi_{i}(\ldots x_{i+m-1} x_{i+m})$ as the $p^i$ coefficient of the $p$-adic representation   of $\phi( \sum_{i=-\infty}^{i+m}x_i p^{i})$.
It is easy to see that $\phi_i$ is well defined and  that  relation (\ref{fgb}) is satisfied. Moreover, for all  $i \in \mathbb{Z},\;
\phi_i$ is   bijective on the last variable.
\hfill $\Box$

 \vspace{1em}

 \noindent {\bf Example}: If $\phi= f_{1/p, 0}$,
 then $\|\phi(x)- \phi(y)\|_p= p \|x- y\|_p \mbox { for all } x, y \mbox { in } \mathbb{Q}_p.$  Furthermore,
 we have
 $$\phi(x)=\sum_{i=-\infty}^{+\infty}\phi_{i}(\ldots x_{i} x_{i+1})p^{i},$$
 where $\phi_{i}(\ldots x_{i} x_{i+1})= x_{i+1}$ for all integer $i$.

\vspace{1em}

{ \bf Proof of Theorem \ref{thb}}:
Assume that $\|a \|_p>1$.  Let $\varepsilon >0$ and $(x_n)_{n \in \mathbb{Z}}$ be a $\varepsilon$-pseudo-orbit.

 Put $$w_n= x_{n+1}- f_{a, b} (x_n), \textrm{  for each  } n \in \mathbb{Z} .$$
We have for all $n \geq 1$,
$$x_n= a^n x_0 + a^{n-1} (w_1+ b)+ a^{n-2} (w_2+b) \cdots+ (w_n + b),$$
and
$$x_{-n}= a^{-n}x_0 - a^{-n+1} (w_{-1}+ b)- a^{-n+2} (w_{-2}+b) \cdots- a^{-1}(w_{-n} + b).$$
Since $\|a \|_p>1$,  the series
$\sum_{i=1}^{+\infty} a^{-i} w_i$ is convergent.

Let
\begin{eqnarray}
\label{dfvnn}
x=  x_0+ \sum_{i=1}^{+\infty} a^{-i} w_i.
\end{eqnarray}
Thus,
$$x_n - f_{a,b}^n (x)= - \sum_{i=1}^{+\infty} a^{-i} w_{n+i} \mbox { for all } n \geq 0$$
 and
$$\|x_n - f_{a, b}^n (x)\|_p \leq \sup  \{ \|w_{k}\|_p,\; k\geq 0\} \leq \varepsilon.$$
On the other hand
for all integer $n >1$,
\begin{eqnarray}
\label{dfvnn2}
x_{-n} - f_{a,b}^{-n} (x)= a^{-n} (x_0 - x) - a^{-n+1} w_{-1}- a^{-n+2} w_{-2} \cdots- a^{-1}w_{-n}.
\end{eqnarray}

By (\ref{dfvnn}) and (\ref{dfvnn2}),
we deduce that
$$\|x_{-n} - f_{a,b}^{-n} (x)\|_p \leq  \sup \{ \|w_{k}\|_p,\; k\in\mathbb{Z} \} \leq \varepsilon.$$
Hence $f_{a,b}$ is shadowing.

The case  $\|a \|_p \leq 1$, left to the reader, can be done by the same manner.

\vspace{1em}

Now, let us prove that $f_{a,b}$ is Lipschitz structurally stable. Since $f_{a, b}$ is topologically conjugate to the map $f_{a, 0}$ by the isometry $\phi: \mathbb{Q}_p \to \mathbb{Q}_p$ defined by $\phi(z)= z + \frac{b}{1- a}. $ We can
 suppose that $b=0$.  Now, assume that $ \|a \|_p>1$ and put $f_{a, 0}= f_a$.
 Let $g= f_a+\psi: \mathbb{Q}_p \to \mathbb{Q}_p$
 such that $\psi: \mathbb{Q}_p \to \mathbb{Q}_p$ is a bounded $\delta$-Lipschitz map, where $\|\psi \|_{\infty} < \delta.$

 \vspace{1em}
 {\bf Claim 1:} The map $g$ is an homeomorphism for $0 < \delta <\|a\|_p - \frac{\|a\|_p}{p}$.

Indeed, let $ x, y \in \mathbb{Q}_p$. Then
$$
\begin{array}{rcl}
(\|a\|_p- \delta) \|x- y\|_p & \leq & \|g(x)- g(y)\|_p \\
                           & \leq & \max\{\|a\|_p \|x- y\|_p,\delta \|x- y\|_p\}.
                           \end{array}$$
                           Since
                           $$0 <\delta < \|a\|_p - \frac{\|a\|_p}{p},$$
                           we deduce that
\begin{eqnarray}
\label{gnj}
\|g(x)- g(y)\|_p = \|a\|_p \|x- y\|_p.
\end{eqnarray}
 Hence, $g$ is continuous and injective.

On the other hand, let $z$ in  $\mathbb{Q}_p$ and
consider the map
$r:  \mathbb{Q}_p \to  \mathbb{Q}_p$ defined by
$$ r(x)= \frac{z}{a}- \frac{\psi(x)}{a}.$$
For all $x, y \in  \mathbb{Q}_p$, we have
 $$ \| r(x)- r(y)\|_p \leq \frac{\delta}{\|a\|_p} \|x- y\|_p.$$
  Hence $r$ is a contraction in the complete metric space
 $\mathbb{Q}_p$.  Thus, $r$  have a fixed point $x$ in $\mathbb{Q}_p$.  Therefore, $g(x)=z$ and we deduce that $g$ is  surjective. This proves  Claim 1.

 \vspace{1em}

 {\bf Case 1}:  $a= 1/p$.

Let
$h(x)=
\sum_{i=-\infty}^{+\infty} h_{i}(x) p^i$ where
$h_i (x) \in \{0,\ldots, p-1\}$ for all $i \in \mathbb{Z}$
and
$$h_i(x)= (g^{i}(x))_0 \mbox { for all } i \in \mathbb{Z}.$$
That is
\begin{eqnarray}
\label{lbcs}
 \ldots  h_{-2}(x) h_{-1}(x). h_0 (x) h_1(x) \ldots =    \ldots (g^{-2}(x))_0 \; (g^{-1}(x))_0. x_0 \; (g(x))_0 \ldots .
 \end{eqnarray}
It is easy to check that
\begin{eqnarray}
\label{fxz}
f_{1/p} \circ h= h \circ g.
\end{eqnarray}

\vspace{0.5em}

{\bf Claim 2:}  $h$ is  an isometry.

\vspace{0.5em}

First, observe that by (\ref{gnj}) we deduce that $g^{-1}: \mathbb{Q}_p \to
\mathbb{Q}_p$ is a contraction.  Hence, there exists a unique fixed point $z^* = \sum_{i=t}^{+\infty} z^*_i p^i $ in $\mathbb{Q}_p,\; t \leq 0$ such that
$g^{-1}(z^*)=z^*$.
{Hence,}
\begin{eqnarray}
\label{gbn}
 z^*+ p \psi(z^*)= p z^*.
 \end{eqnarray}
Since $\|\psi \|_{\infty} \leq \delta < 1$, we deduce that
$z^*_i= (z^*+ p \psi(z^*))_i$ for all $i \leq 0$.  Hence,
by (\ref{gbn}), we have
$$z^*_i=  (pz^*)_i \mbox { for all } i \leq 0.$$
Thus
\begin{eqnarray}
\label{lpcx}
z^*_0= z^*_{-1}= \ldots= z^*_{t}= z^*_{t}=0.
\end{eqnarray}

Now,  if $x \in \mathbb{Q}_p$, then
by (\ref {gnj}), we deduce that
$$\|g^{-n}(x) - z^*\|_p= p^{-n}\|x - z^*\|_p \mbox { for all } n \geq 0.$$
 Thus, $$\lim_{n \to \infty} g^{-n}(x) = z^*.$$
We deduce   that there exists a positive integer $n_0$  such that
$$(g^{-n}(x))_0= z^*_0 = 0 \mbox { for all } n > n_0.$$
Hence $h(x) $ belongs to $\mathbb{Q}_p$ and $h$ is well defined.

\vspace{0.5em}

Now, let $x = \sum_{i=-\infty}^{+\infty} x_i p^i$ and  $ y  = \sum_{i=-\infty}^{+\infty} y_i p^i$ be two elements of
$\mathbb{Q}_p$ such that $\|x- y\|_p = p^{l},\; l \in \mathbb{Z}$.
By (\ref {gnj}), we deduce that for all integer $n \in \mathbb{Z}$,
\begin{eqnarray}
\|g^n (x)- g^n(y)\|_p = p^{n+l}.
\end{eqnarray}
In particular,
$$\|g^{-l} (x)- g^{-l}(y)\|_p = 1 \mbox { and } \|g^n (x)- g^n(y)\|_p < 1 \mbox { for all } n <-l  .$$
 Hence,
 $$(g^{-l} (x))_{0} \neq (g^{-l} (y))_{0} \mbox { and } (g^{n} (x))_{0} = (g^{n} (y))_{0} \mbox { for all } n <-l .$$

 Thus, $\|h(x) - h(y)\|_p= p^{l}=\|x - y\|_p $ and we obtain the claim.

 \vspace{1em}

 Using Claim 2, we deduce that $h$ is an injective continuous map.

 \vspace{1em}

To prove that $h$ is an homeomorphism, it suffices to prove that $h$ is surjective.

Indeed, let
 $y= \sum_{i=l}^{+\infty} y_i p^i \in \mathbb{Q}_p$ where $y_l >0$.
 We aim to find
 $$x= \sum_{i= s}^{+\infty} x_i p^i \in \mathbb{Q}_p \mbox { such that } h(x)= y,$$
 where $x_s >0.$
 Since $h$ is an isometry, we have
 $$\|h(x)- h(z^{*})\|_p= \|x - z^*\|_p.$$
 By (\ref{fxz}), we deduce that
 $\frac{1}{p} h(z^*)=  h(z^*)$.  Hence,
 $h(z^*)=0$.
  Therefore,
 \begin{eqnarray}
 \label{lza}
 \|h(x)\|_p= \|x - z^*\|_p.
 \end{eqnarray}
 Using (\ref{lpcx}) and (\ref{lza}), we deduce that
  $s=l$. Hence,
 $x= \sum_{i= l}^{+\infty} x_i p^i,$ where $x_l >0$.
 By (\ref{lbcs}), we deduce that

 \begin{eqnarray}
 \label{jkn}
 y_i = (g^i(x))_0 \mbox { for all } i \geq l.
 \end{eqnarray}
 Hence, $x_0= y_0$.

Assume that $l \leq 0$. By (\ref{gnj}), we deduce that
 for all  $x, y \in \mathbb{Q}_p,$
 $$\|g^l(x)- g^{l} (y) \|_p =\|a\|_p^{-l} \|x- y \|_p=
 p^{l} \|x- y \|_p.
 $$
 By Lemma \ref{swvb} applied to $\phi= g^l$, we have
 $$y_l = (g^l(x))_0= \phi_0 (     \ldots ,0,0,x_l).$$
 Using the fact that $\phi_0$ is bijective in the last variable, we find $x_l$.
Using (\ref{jkn}), we obtain recursively
 $$x_{i},\; i >l.$$
 Hence, $h$ is surjective. Since $h$ is an isometry,  it follows that $h^{-1}$ is  continuous and hence $f_{1/p, 0} $  is Lipschitz structurally stable.

 Observe that the only thing used about $g$ is that
 $\|g(x)- g(y)\|_p= p \|x- y\|_p$ for all $x, y \in \mathbb{Q}_p$.

 Hence, we deduce that if $\|a\|_p= p$, then $f_{a,0}$ is topologically conjugate to $f_{1/p, 0}$ and the conjugation is an isometry. Thus, $f_{a,0}$ is Lipschitz structurally stable.

 \vspace{1em}

 {\bf Case 2}: $a= 1/ p^k$ where $k \geq 2$ is an integer.

 This case can be done by the same way
by taking  the map $h$ defined by
$$ \ldots  h_{-1}(x).h_0(x) \ldots= \ldots (g^{-1} (x))(0, k-1) \;. \;x(0, k-1) \; (g (x))(0, k-1)  \ldots $$
where
$$(g^{i} (x))(0, k-1)= (g^{i} (x))_0 \ldots (g^{i} (x))_{k-1} \mbox { for all } i \in \mathbb{Z}.$$

The case  $\|a \|_p \leq 1$, left to the reader, can be done by the same manner. \hfill $\Box$

\begin{theorem}
\label{mzx}
Let   $f : \mathbb{Q}_p \to \mathbb{Q}_p$ be an homeomorphism. If $f$  is a contraction (Lipschitz map with Lipschitz constant $\leq p^{-1}$) or a dilatation ($f^{-1}$ a contraction),
then $f$ is shadowing and Lipschitz structurally stable.
 \end{theorem}

 \begin{remark}
 The last result is more general  than Theorem \ref{thb}. However, in the  proof, we will use the fixed point theorem to find the real orbit that shadow the $\delta$-pseudo-orbit and also to prove that the conjugation is bijective. In the opposite of the case $f_{a, b}$, here we don't exhibit explicitly the conjugation and we don't know if it is an isometry.
 \end{remark}

 {\bf Proof of Theorem \ref{mzx}.}
 The proof will be done in the case where $f$ is a dilatation and is  analogous for the other case.

Let $k \geq 1$ be an integer and assume that
 \begin{eqnarray}
\label{svz}
\|f(x)- f(y)\|_p \geq  p^k \|x- y\|_p \mbox { for all } x, y \mbox { in } \mathbb{Q}_p.
\end{eqnarray}
 Let $g= f+\psi: \mathbb{Q}_p \to \mathbb{Q}_p$
 such that $\psi: \mathbb{Q}_p \to \mathbb{Q}_p$ is a bounded $\delta$-Lipschitz map, where $\|\psi \|_{\infty} < \delta < p^{k}- p^{k-1}$.

 \vspace{0.5em}

 {\bf Claim 1:} The map $g$ is a shadowing homemorphism.

 \vspace{1em}

 Indeed,  we deduce using triangular inequality that
 \begin{eqnarray}
 \label{zab}
 \|g(x)- g(y)\|_p \geq  p^k \|x- y\|_p \mbox { for all } x, y  \in \mathbb{Q}_p.
 \end{eqnarray}
 Hence, $g$ is injective.

 On the other hand, given $z$ in  $\mathbb{Q}_p$ and applying the fixed point theorem to the  map
$s:  \mathbb{Q}_p \to  \mathbb{Q}_p$ (which is a contraction), defined by
$$ s(x)= f^{-1} (z - \psi(x)),$$
 we deduce that $g$ is surjective and hence $g$ is an homeomorphism

\vspace{1em}

 Let $E= \mathbb{Q}_p^{\mathbb{N}}$ be the set of sequences $(x_{n})_{n \geq 0}$ such that $x_{n}  \in \mathbb{Q}_p$ for all integer $n \geq 0$.
  We endow $E$ with the metric $d$ defined by:

 $$d((x_{n})_{n \geq 0}, (y_{n})_{n \geq 0})= \sup \{\|x_n- y_n \|_p,\; n\in\mathbb{N}\}.$$
It is not difficult to prove that $(E, d)$ is a complete metric space.

 Now, let $\varepsilon>0$ and consider a $\varepsilon$-pseudo-orbit $x= (x_n)_{n \in \mathbb{Z}}$ and define the function
 $\Phi: E \to E$ by
 $$\Phi ((y_n)_{n \geq 0})= (g^{-1}(x_{n+1}+ y_{n+1})- x_n)_{n \geq 0}.$$
 By (\ref{zab}), we deduce that
 the map $\Phi$ is a contraction.
 Moreover, if $d((y_n)_{n \geq 0}, 0) \leq \varepsilon$,
 then $$\|g^{-1}(x_{n+1}+ y_{n+1})- x_n\|_p  \leq\max \left\{ \|g^{-1}(x_{n+1}+ y_{n+1})- g^{-1}(x_{n+1})\|_p, \|g^{-1}(x_{n+1})- x_n\|_p \right\}.$$

 Hence,
 $$\|g^{-1}(x_{n+1}+ y_{n+1})- x_n\|_p  \leq
\max \left\{p^{-k} \|y_{n+1}\|_p, p^{-k}\|x_{n+1}- g(x_n)\|_p\right\}  \leq p^{-k} \varepsilon.$$

Hence, $\Phi$ sends the closed ball of $E$ of center $0$ and radius $\varepsilon$ into the same ball. Thus, $\phi$ has a fixed point $(y^*_n)_{n \geq 0}$ in this ball.
 We deduce that
 $g(x_n+ y^*_n)= x_{n+1}+ y^*_{n+1} \mbox { for all } n \geq 0.$
 Hence, $ g^n (x_0 + y^*_0)- x_n= y^*_n \mbox { for all } n \geq 0.$
 Thus,
\begin{eqnarray}
\label{szqu}
\|g^n (x_0 + y^*_0)- x_n \|_p \leq \varepsilon  \mbox { for all } n \geq 0.
\end{eqnarray}

On the other hand,  we deduce by (\ref{zab}) that $g^{-1}$ is a contraction.
Hence, for all $n \geq 0$
$$\|g^{-1} (x_{n+1})- x_n \|_p \leq \frac{\varepsilon}{p^k} \leq \varepsilon.$$
Thus, the sequence $(x_{-n})_{n \geq 0}$ is a $\varepsilon$-pseudo-orbit for $g^{-1}$.
We deduce as done in Proposition \ref{lipc} that $g^{-1}$ is shadowing and
 $\|g^{-n} (x_0) - x_{-n} \|_p \leq \varepsilon.$
Thus, for all integer $n \geq 0$,
\begin{eqnarray}
\label{mxaqd} \|g^{-n} (x_0 + y^*_0 ) - x_{-n} \|_p \leq
 \max \left\{ \frac{\varepsilon}{p^{kn}}, \varepsilon \right\} \leq \varepsilon.
 \end{eqnarray}

By (\ref{szqu}) and (\ref{mxaqd}), we deduce that
 $g$ is shadowing and we obtain the claim and also the fact that $f$ is shadowing.

\vspace{0.5em}

{\bf Claim 2:} The map $f$ is Lipschitz structurally stable.

The proof of Claim 2 is similar to the proof of Lemma \ref{str}.

For $x\in\mathbb{Q}_p$, we have that $$\| g^{n+1}(x)-f(g^n(x))\|_p=\| (g-f)(g^n(x))\|_p < \delta.$$
Since $f$ is shadowing, there exists  $y= h(x) \in\mathbb{Z}_p$ such that
for all $n\in\mathbb{\Z}$,
\begin{eqnarray}
\label{f1}
\| g^n(x)-f^n(h(x))\|_p \leq \delta.
\end{eqnarray}

Since $f$ is a dilatation, we deduce that
$h(x)$ is unique and $f\circ h=h\circ g$.

Now, since $g$ is a dilatation, by Claim 1, $g$ is shadowing, we deduce by changing the roles of $f$ and $g$ that $h$ is bijective.

In order to prove that $h$ and $h^{-1} $ are  continuous,
let $x \in \mathbb{Q}_p$. Since $g$ is continuous, given an integer $N >0$, there exists a real number $\alpha >0$ such that for all $y \in \mathbb{Q}_p$ if $\|x-y\|_p \leq \alpha$, then
\begin{eqnarray}
\label{dfbl}
\|g^n(x)- g^n(y) \|_p \leq \delta, \textrm{ for all } 0 \leq n \leq  N.
\end{eqnarray}
Hence, if $\|x-y\|_p \leq \alpha$, then we deduce
by (\ref{f1}) and (\ref{dfbl}) and the fact that $f$ is a dilatation  that
\begin{eqnarray}
\|f^n(h(x))- f^n(h(y)) \|_p \leq \delta \mbox { for all }  0 \leq n \leq  N.
\end{eqnarray}
Thus, $$\|h(x)- h(y) \|_p < \frac{\delta}{p^{k N}}$$ and we deduce that $h$ is continuous.
Changing the roles of $f$ and $g$, we deduce that $h^{-1}$ is continuous and we are done
\hfill $\Box$

\vspace{1 em}

{\bf Questions:}   Let $f : \mathbb{Q}_p \to \mathbb{Q}_p$ be an homeomorphism.
\begin{enumerate}
\item Assuming that  $f$ is $1$-Lipschitz,
can $f$ be shadowing or Lipschitz structurally stable?
\item Suppose that  $f$ is a shadowing and expansive map. Is $f$
Lipschitz structurally stable?
\end{enumerate}

\vspace{1 em}

\noindent {\bf Acknowledgment}: The third author would like to thank Fabien Durand for fruitful discussions.
The first author was partially supported by Fapesp project 2013/24541-0.
The second author was supported by Fapesp project 2015/26161-6 and Capes project 88887.468127/2019-00.
The third author was partially supported by CNPq project \\ 311018/2018-1  and by Fapesp project  2019/10269-3.

\end{document}